\documentclass[twoside, 10 pt]{article}

\usepackage{amsfonts}
\usepackage{amsmath}
\usepackage{amscd}
\usepackage{amssymb}
\usepackage{latexsym}
\usepackage{graphicx}

\newtheorem{Th}{THEOREM}[section]
\newtheorem{lem}[Th]{LEMMA}

\newtheorem{df}[Th]{DEFINITION}
\newtheorem{prop}[Th]{PROPOSITION}

\newcommand{\Rup}{{\mathbb R}^5_{>0}}
\newcommand{\Psl}{PSL(2,\H)}
\newcommand{\Mr}{{\cal M}_r}
\newcommand{\Pn}{{\cal P}_r}
\newcommand{\Lr}{{\cal L}_r}
\newcommand{\Ur}{{\cal U}_r}
\newcommand{\HP}{{\mathbb H}{\mathbb P}}
\newcommand{\CP}{{\mathbb C}{\mathbb P}}
\newcommand{\3}{\Mr^{(3)}}
\newcommand{\dl}{\delta}
\newcommand{\e}{\epsilon}
\newcommand{\vt}{\vartheta}
\newcommand{\spn}{\langle}
\newcommand{\es}{\rangle}

\newcommand{\sq}{\hfill{\Box}\medskip}
\newcommand{\pf}{{\noindent \it Proof.{\mbox{  }}}}

\def \RP{{\mathbb R\mathbb P}}
\def \C{{\mathbb C}}
\def \E{{\mathbb E}}
\def \H{{\mathbb H}}

\def \R{{\mathbb R}}

\def \L{{\cal L}}
\def \M{{\cal M}}

\def \S{{\Sigma}}

\def \P{{\mathbb P}}

\def \e{{\varepsilon}}
\def \g{{\gamma}}

\def \l{{\lambda}}

\def \su{{\mathfrak s}{\mathfrak u}}
\def \calO{{\cal O}}
\def \calH{{\cal H}}

\def \s{{\sigma}}
\def \t{{\theta}}

\def \HP{{\mathbb H}{\mathbb P}}
\def \dim{{\rm dim}}
\def \cprime{$'$}

\long\def\startmy#1\endmy{}

\date{}

\begin{document}

\title {The Geometry of Polygons in $\R^5$ and Quaternions} 

\author {Philip Foth and Guadalupe Lozano}

\date {\today}
\maketitle 

\begin{abstract}

We consider the moduli space ${\cal M}_r$ of polygons with fixed side lengths in 
five-dimensional Euclidean space.  We analyze the local structure of its 
singularities and exhibit a real-analytic equivalence between ${\cal M}_r$ 
and a weighted quotient of $n$-fold products of the quaternionic projective 
line $\H\P^1$ by the diagonal $PSL(2,\H)$-action. We explore the relation 
between ${\cal M}_r$ and the fixed point set of an anti-symplectic involution 
on a GIT quotient $Gr_{\C}(2,4)^n/SL(4,\C)$. We generalize the Gel{\cprime}fand-MacPherson 
correspondence to more general complex
Grassmannians and to the quaternionic context, and realize our space ${\cal M}_r$ 
as a quotient of a subspace in the quaternionic Grassmannian $Gr_{\mathbb H}(2,n)$ 
by the action of the group $Sp(1)^n$. We also give analogues of the Gel{\cprime}fand-Tsetlin 
coordinates on the space of quaternionic Hermitean matrices and briefly describe generalized 
action-angle coordinates on ${\cal M}_r$.  

\end{abstract}

\section{Introduction}

\footnotetext{{\it AMS subj. class.}: 
\ \ primary 58D29, secondary 53D20, 14D20.}
\footnotetext{{\it Keywords:} polygons, quaternions, invariants, 
reduction, lagrangian, Grassmannians} 

Spaces of polygons and linkages provide a useful class of geometric examples, as their structure can be effectively visualized using linkages themselves. For example, the spaces of polygons ${\cal M}_r^{(3)}$ with fixed side lengths in the Euclidean space ${\mathbb E}^3$ (studied by Klyachko \cite{Kly}, Kapovich and Millson \cite{KM}, Hausmann and Knutson \cite{HK} among others), carry interesting symplectic and K\"ahler structures.  These spaces also possess natural action-angle coordinates, which can be obtained using symplectic quotients of complex Grassmannians by the action of maximal tori.  Moreover, for a rational choice of side lengths, these spaces can be identified with the Mumford quotients of $({\CP^1})^n$ by the diagonal action of $SL(2, \C)$. By a theorem of Kapranov \cite{Kapr}, the Deligne-Mumford compactification $\overline{M}_{0,n}$ of the moduli space of $n$-marked projective lines dominates all of ${\cal M}_r^{(3)}$. This map was explicitly constructed in \cite{F2}. Hu in \cite{Hu} used polygon spaces to study the K\"ahler cone of the space $\overline{M}_{0,n}$. He also constructed an explicit sequence of blow-ups leading from ${\cal M}_r^{(3)}$ to $\overline{M}_{0,n}$.

	In the present paper we study the moduli spaces $\M_r$ of polygons in ${\mathbb E}^5$ and relate them to certain quotients of products of quaternionic projective lines by the diagonal action of $SL(2, {\mathbb H})$. 

	As in the case of $\3$, each $\Mr$ can be naturally associated to the $r$-weighted quotient of $(S^4)^n$ by the diagonal action of $SO(5)$, where $r\in R^n_{>0}$ is the weight vector prescribing the side lengths.   

	As constructed, the spaces $\Mr$ turn out to be singular, with the singular locus isomorphic to a quotient of ${\cal M}_r^{(3)}$ by an involution.  In Section 2, we study the local structures of these singularities and show that they come in two types. A neighborhood of a planar $n$-gon is analytically isomorphic to $\R^{n-3}\times [(\R^3)^{n-3}/SO(3)]$ and a neighborhood of a spatial $n$-gon is isomorphic to $\R^{2n-6}\times [(\R^2)^{n-4}/SO(2)]$.  In both cases, the first factor corresponds to deformations along the smooth singular locus.

	In Section 3 we establish a real-analytic isomorphism between 
$\Mr$ and the weighted quotient of $n$-copies of the quaternionic projective line, 
$(\HP^1)^n$, by $PSL(2,\H)$ for a weight vector $r$.  We use stable measures
due to Leeb and Millson \cite{LeebM} and apply them to $S^4$, the geometric
boundary of the symmetric space $B^5\simeq SL(2, {\mathbb H})/Sp(2)$. We also use 
properties of the conformal barycenter described by Douady and Earle in \cite{DE}.

	The main result in Section 4 generalizes the Gel{\cprime}fand-MacPherson correspondence to the quaternionic context and shows that our spaces ${\cal M}_r$ may be realized as certain quotients of the quaternionic Grassmannian $Gr_{\mathbb H}(2,n)$ by the action of the group $Sp(1)^n$. These correspondence together with the analogue of the Gel{\cprime}fand-Tsetlin coordinates on the space of quaternionic Hermitean matrices, described in Section 5, allows us to use the results of \cite{HK} and introduce a certain generalization of action-angle coordinates on ${\cal M}_r$, as shown in Section 7.

	One of the major results of the paper appears in Section 6.  There we relate the space ${\cal M}_r$ with a fixed point set of an involution on the Mumford quotient $Gr_{\C}(2,4)^n/SL(4,\C)$ for an appropriate choice of linearization.  The involution in question arises both from multiplication by ${\bf j}$ on the space $\C^{2m}$ identified with ${\mathbb H}^m$, and the involution of $SL(4,\C)$ defining the real form $SL(2, {\mathbb H})$.  In this section we also extend the Gel{\cprime}fand-MacPherson correspondence for projective spaces studied in \cite{Kapr} to more general complex Grassmannians.

	As we show, a certain number of methods for studying ${\cal M}_r$ can essentially be extended from those for the polygon spaces studied previously. However, we always keep our focus on a number of new phenomena, pertinent to our specific objects and settings.
\bigskip

\noindent{\bf Acknowledgments.} We would like to thank Sam Evens, Yi Hu, Jiang-Hua Lu, 
John Millson, and Reyer Sjamaar for useful conversations and correspondence. The first 
author is partially supported by NSF grant DMS-0072520.
\bigskip

\section {Real analytic structure of $\Mr$}

	Our goal in this section is to prove two main results related to the singularities of $\Mr$.  First, we will show that, for generically chosen $r$, the singular locus, $D_r$, of $\Mr$ is a real orbifold of dimension $2n-6$.  It is naturally isomorphic to $\3/({\mathbb Z}/{2{\mathbb Z}})$, where $\3$ denotes the moduli space of closed polygons in ${\mathbb E}^3$.  We further show that, for $n\geq 4$, $\Mr$ possesses two types of singularities, each one of which is locally equivalent to a neighborhood of zero in $\R^{2n-6}\times [(\R^2)^{n-4}/SO(2)]$ or in $\R^{n-3}\times [(\R^3)^{n-3}/SO(3)]$.  We conclude the section with a geometric description of singularities for $n=5,6$.	
	
	Fix an ordered $n$-tuple $r=(r_1,...,r_n)$ of positive real numbers and let $\Lr$ denote the collection of $n$-sided polygonal linkages in ${\mathbb E}^5$ with fixed side lengths $r_i$, $i=1,...,n$.  Clearly, a polygonal linkage $P$ lying in $\Lr$ determines an $n$-tuple of unit vectors $U=(u_1,...,u_n)$ in $(S^4)^n$.  Conversely, such an $n$-tuple uniquely determines a linkage $P$ with side lengths prescribed by the $r_i$.  Hence, by considering the $SO(5)$-equivariant map:
$$\phi:(S^4)^n\rightarrow {\mathbb R}^5, \ U\mapsto \sum_{i=1}^n {r_iu_i}$$
we can realize the space $\Pn$ of $5$-dimensional (closed) $n$-gons as the zero locus, $\Ur:=\phi^{-1}(0)$ of the above map and obtain $\Mr$ as the quotient of $\Ur$ by the diagonal action of $SO(5)$ on $\Ur$:
$$\Mr=\Ur/SO(5).$$

\begin{df}
A polygon $P\in \Pn$ (or its representative, $U\in\Ur$) is called degenerate if it is stabilized by a non-trivial subgroup of $SO(5)$.
\end{df}

	We first note that if $P$ is a non-degenerate $n$-gon then it has trivial stabilizer in $SO(5)$ or, equivalently, any collection of unit vectors $U=(u_1,...,u_n)$ in $(S^4)^n$ representing $P$ spans at least ${\mathbb R}^4$.  It follows that any infinitesimal deformation of $P$ in $\Pn$ yields another non-degenerate polygon, $P_{\dl}$, in $\Pn$.  Of course, $P_{\dl}$ and $P$ will be equivalent if and only if they can be joined by an integral curve of some fundamental vector field of the $SO(5)$-action on $\Pn$. Thus, any non-degenerate polygon has a sufficiently  small neighborhood in $\Pn$ in which $SO(5)$ acts freely.  Consequently, $\Mr$ is locally smooth at $[P]$ as long as $P$ is non-degenerate.

	A weight vector $r$ is called $admissible$ if the corresponding space $\Pn$ is non-empty or, equivalently, $2r_j\leq \sum r_i$.

	Unless noted otherwise, we will restrict our attention to those admissible $r$ which also satisfy the trivial non-degeneracy condition, namely $r_1\pm r_2\pm\cdots\pm r_n\neq 0$.  This last condition eliminates degenerate polygons contained in a straight line, that is, all those stabilized by a subgroup of $SO(5)$ isomorphic to $SO(4)$.  It is easy, however, to see that such trivial polygons do not exhaust the collection of degenerate ones.  For example, if $P\in \Pn$ is such that its edges merely span a linear subspace equivalent to $\R^3$ or $\R^2$, then $P$ is stabilized by $H\subseteq SO(5)$ where $H\simeq SO(2)$ or $H\simeq SO(3)$, respectively and thus is degenerate.  In turns out that these two are, in fact, the only types of degeneracies allowed by our choice of $r$.

\begin{df}
A degenerate polygon is of type $k$, if it is stabilized by $H_k\subseteq SO(5)$ where $H_k\simeq SO(k)$.  Accordingly, a type $k$ singularity in $\Mr$ is one which lifts to a degenerate polygon of type $k$ in $\Ur$.
\end{df}

\begin{lem}
The only degenerate $n$-gons in $\Ur$ are of type $2$ and $3$.
\end{lem}

\pf  If $U=(u_1,...,u_n)$ is stabilized by $H$, then, perhaps after a change of basis, any element of $H$ may be represented by a block-diagonal matrix $h$ having a $k\times k$ identity block, where $k=\dim\spn u_1,...,u_n\es$.  The fact that $h$ must lie in $SO(5)$ forces the second $(n-k)\times (n-k)$ block to lie in $SO(n-k)$. $\sq$

\startmy 
\hrulefill
EXCISED.  (Extended proof of Lemma) Suppose a polygon $P$ 
is degenerate.  Then its edges fail to span ${\mathbb R^5}$.  
Let $A$ be the subspace of ${\mathbb R^5}$ spanned by the 
edges of $P$ and let $\{a_1,..., a_k\}$ be an  basis for $A$.  
(Note that $k$ is necessarily less than $5$). Then $B\in SO(5)$ 
stabilizes $P$ if and only if $B|_A=id$.  So, with respect to 
the basis $\{a_1,....,a_k,b_1,...b_{5-k}\}$ of ${\mathbb R^5}$, 
the first $k$ columns of $B$ looks like $(0,...0,1,0,...0)^T$, 
where the $1$is in the $k^{\mbox{{\small th}}}$ spot.  Since 
$B$ satisfies $BB^T=I$, then the first $k$ rows of $B$ are the 
transposes of the first $k$ columns.  This says that $B$ is 
block diagonal with the top-left block being just a $k\times k$ 
identity matrix. The second $n-k\times n-k$ block, call it $D$, 
must also satisfy $DD^T=I$ since this is necessary for $B$ to 
satisfy this same orthogonality condition.  The above 
orthogonality condition implies $\det(D)=\pm 1$.  We know $\det(B)=1$ so $D$ must either belong to $SO(5)$ or to $O(5)$ at worst, depending on the value of $k$.  If $k=1$, $\det(B)=1(\det(D))$ which implies  $\det(D)=1$ and hence $D\in SO(4)$ (this is the trivial degeneracy case).  If $k=2$, $\det(B)=(\det(upper block=2\times 2 id)(\det(D))=1(\det(D))$ which implies  $\det(D)=1$ and hence $D\in SO(3)$.  Since the upper block is always and identity matrix, and the total determinant is the product of the determinant of the blocks, we see that $D$ must alway have determinant 1, hence  lie in $SO(n-k)$.// 
$$ B=
\left( \begin{array}{ccccc} 
         1 & 0 & 0 & 0 & 0\\
         0 & 1 & 0 & 0 & 0\\
         0 & 0 & 1 & 0 & 0\\
	   0 & 0 & 0 & a & c\\
	   0 & 0 & 0 & c & b
\end{array}\right)\ , ab-c^2=1
$$
\hrulefill
\endmy

	Recall that $\3$ denotes the moduli space of closed polygons in $\E^3$.

\begin{prop}
The singular points in $\Mr$ form a $2n-6$ dimensional orbifold, $D_r$, which is naturally isomorphic to $\3/({\mathbb Z}/{2{\mathbb Z}})$.
\end{prop}

\pf  Define an involution $\vt$ in $\3$ by $\vt([P])=[\s(P)]$, where $P$ is any representative of $[P]\in \3$ and $\s$ is an arbitrary reflection about a plane in $\R^3$.  Note that $\vt$ is well defined, as for any two reflections $\s_1$, $\s_2$ in $\R^3$ $[\s_1(P)]=[\s_2(P)]$.  Also $[\s_1(P)]=[\s_1(Q)]$, if $[P]=[Q]$.  So $\vt$ defines a ${\mathbb Z}/{2{\mathbb Z}}$-action on $\3$ fixing all planar polygons.  Now suppose $[P]$ and $[Q]$ share a ${\mathbb Z}/{2{\mathbb Z}}$-orbit.  If $[P]$ and $[Q]$ are planar, then $[P]=[Q]$, that is, $P$ and $Q$ lie in the same $SO(5)$-orbit.  Otherwise, assume that $P$ and $Q$ have been chosen so that they span the same $3$-dimensional subspace $T\in \R^5$.  As $P$ and $Q$ differ by a reflection in $T$, we may choose an orthonormal basis $\{e_1, e_2\}$ for the plane fixed by the reflection and extend it to a basis of $T$ by choosing $e_3$ to be normal to the fixed plane.  Then, in any basis $\{e_1, e_2, e_3, \cdot, \cdot\}$ of $\R^5$, the matrix
$$ B=
\left( \begin{array}{ccccc} 
         1 & 0 & 0 & 0 & 0\\
         0 & 1 & 0 & 0 & 0\\
         0 & 0 & -1 & 0 & 0\\
	   0 & 0 & 0 & a & c\\
	   0 & 0 & 0 & c & b
\end{array}\right)\ , ab-c^2=-1
$$
is an element of $SO(5)$ mapping $P$ to $Q$. $\sq$

\startmy 
\hrulefill EXCISED.

Shorter pf: 
It is now straight-forward to check that if $[P]$ and $[Q]$ lie in the same ${\mathbb Z}/{2{\mathbb Z}}$-orbit, corresponding $P$ and $Q$ lie in the same $SO(5)$-orbit.

\hrulefill
\endmy

 	At this point, we note that for $n=4$, $\Mr=D_r$, as every $4$-sided polygon must lie entirely within a $3$-dimensional subspace in $\R^5$.  Since the moduli space of $4$-gons in $\R^3$ is homeomorphic to a $2$-sphere (for $r$ satisfying the trivial non-degeneracy condition), it follows that in this case $\Mr\simeq S^2/({\mathbb Z}/{2{\mathbb Z}})$; that is, $\Mr$ is isomorphic to a closed disk.

	We also remark that Kamiyama in \cite{Kami2} showed that for equilateral hexagons, the space $\Mr$ is homeomorphic to $S^9$.  In \cite{Kami1} he computed the Euler characteristic of the space of equilateral septagons, which turns out to be equal to $-7$. Also, Schoenberg in \cite{Schoen} showed that, for the case of equilateral septagons, the quotient of $\Mr$ by a natural involution (defined as in Proposition 2.4 above) is homeomorphic to $S^{13}$.

	We now turn to the analysis of the local structure of $\Mr$ near singular points.  Our strategy makes use of slices.  These are essentially generalized cross-sections which allow us to understand the structure of $\Mr$ near its singularities by studying the action of point stabilizers along directions transverse to local orbits. 

	Recall that for a compact Lie group $G$, a $G$-space $M$ and $q\in M$, a $slice$, $S$, around $q$ is a subset of $M$ containing $q$ and satisfying the following conditions:
\begin{itemize}
\item[{\bf a.}] $S$ is closed in $G\cdot S$;
\item[{\bf b.}] $G\cdot S$ is an open neighborhood of the orbit of $q$, $G\cdot q$;
\item[{\bf c.}] $G_q\cdot S=S$, ($G_q$ denotes the stabilizer of $q$ in $G$);
\item[{\bf d.}] $gS\cap S \neq \varnothing \Rightarrow g \in G_q$.
\end{itemize}
	We will begin by considering the type of singularity arising from degenerate polygons of type 2.  The case of type 3 degeneracies is analogous and will be outlined at the end of this section. 
 
	Let $P_0$ denote a degenerate polygon of type 2.  Since we are interested in a local result, we may, if necessary, permute the edges of $P_0$ so that the first three span some 3-dimensional subspace of $\R^3$.  Assume next that $U_0=(u_1^0,...,u_n^0)$, $u_i^0 \in S^4$ represents $P_0$, so that $\dim(\spn u_1,u_2,u_3\es)=3$.  Because $[U_0]$ is defined up to the action of $SO(5)$, we may assume that $u_1^0=e_1$, $u_2^0\in \spn e_1,e_2\es^+$ and $u_3^0\in \spn e_1,e_2,e_3\es^+$, where $e_i, \ i=1,...,5$ denote the standard basis vectors of $\R^5$.  (Here, the $+$ sign selects vectors with positive coordinates in the $e_2, \  e_3$ directions.)

	Define,
\begin{eqnarray}
S=\{U=(u_1,...,u_n)&\in& (S^4)^n:u_1=e_1, u_2\in \spn e_1,e_2\es, \nonumber \\
                u_3&\in& \spn e_1,e_2, e_3\es, \ U \mbox{ spans }\spn e_1,e_2, e_3\es \}.\nonumber
\end{eqnarray}

	Then 
$$S_0=\{U\in S: r\cdot U=\sum_{i=1}^n{r_i u_i}=0\}$$
is a smooth $4n-14$ dimensional submanifold of $\Ur$ and a slice through $U_0$ for the $SO(5)$-action on $\Ur$.  Note that the requirement that $U\in S$ spans $\R^3$ ensures that no degenerate polygons of type $3$ lie in $S$.  

\startmy
\hrulefill
Intuitively, a slice should contain at least one representative from each of the $SO(5)$-orbits near $u_0$ that is, we must make sure that each and all infinitesimal deformations of $U_0$ are represented in the slice.  Observe that since $U_0$ spans $\R^3$ there exists a neighborhood of $U_0$ transversal to the local orbits (so just comprised of ``non-equivalent" deformations) which does not contain any planar configurations.  As mentioned above, this allows us to choose $S$ containing only type $2$ degeneracies. 
\hrulefill
\endmy

	The fact that $S_0$ is a slice implies that the natural map
$$\varphi: SO(5)\times_{H_2} S_0\rightarrow\Ur, \ [g,U]\mapsto gU$$
is a tube around the $SO(5)$-orbit of $U_0$ that is, $\varphi$ is an $SO(5)$-equivariant diffeomorphism onto a neighborhood of the orbit of $U_0$ in $\Ur$, where the $SO(5)$-action on $SO(5)\times_{H_2} S_0$ is just left translation 
$$g([g',U])=[gg',U], \  \forall g\in SO(5).$$ 
It follows that,
$$\Ur/SO(5)\simeq(SO(5)\times_{H_2} S_0)/SO(5)\simeq S_0/H_2$$
near $U_0$.

\startmy
\hrulefill
Since left translation is a transitive action, the orbit of $[g,U]$ in $SO(5)\times_{H_2} S_0$ is $[SO(5),U]$.  Now, $[SO(5),U]=[SO(5),V]$ if and only if there exist $h \in H_2$ such that $SO(5)h^{-1}=SO(5)$ (any $h$ will do), and $hU=V$, i.e., $U$ and $V$ share the same $H_2$-orbit in $S_0$.  
\hrulefill
\endmy

\startmy
\hrulefill
$SO(5)\times_{H_2} S_0$ denotes the {\it twisted product} of $SO(5)$ and $S_0$, i.e., the orbit space of the $H_2$-action on $SO(5)\times S_0$ given by $h(g,u)=(g h^{-1},hu)$.  We can think of this space as a fiber bundle over $SO(5)/H_2$ with fiber $S_0$, hence the naturality of the map $\varphi$. In fact, a slice can also be characterized as $G_q$-invariant subspace of $M$ for which the map $\varphi$ is a tube around $G\cdot q$ (an equivariant embedding onto some open neighborhood of this orbit).
\hrulefill
\endmy

	Now, as $S$ is clearly smooth near $U_0$ and $0$ is a regular value of the map $\phi|_S$, it follows that $S_0=\phi|_S^{-1}(0)$ is also locally smooth near $U_0$.  Thus, $H_2$ induces a linear action on $T_{U_0}S_0$ via the isotropy representation.

	Let $N(U_0)$ denote a sufficiently small neighborhood of $U_0$ in $S_0$.  In the presence of an $H_2$-equivariant map $\psi:N(U_0)\rightarrow T_{U_0}S_0$, we have that $N(U_0)/H_2\simeq T_{U_0}S_0/H_2$.  Thus, understanding the local structure near the singularity $[P_0]\in \Mr$ amounts to understanding the linear action of $H_2\simeq SO(2)$ on $T_{U_0}S_0$.

\begin{prop}
If $[P_0]$ is a type $2$ singularity, then there exists a neighborhood of $[P_0]$ in $\Mr$ isomorphic to a neighborhood of $0$ in 
$$\R^{2n-6} \times [(\R^2)^{n-4}/SO(2)].$$  
The factor $\R^{2n-6}$ corresponds to infinitesimal deformations along the singular locus, that is, those spanning three dimensions.
\end{prop}

\pf  The vector space $T_{U_0}S_0$ consists of vectors $\e=(\e_1,...,\e_n)$ in $(\R^5)^n$ satisfying the following conditions:
\begin{itemize}
\item[({\bf i})] $\e_i\cdot u_i^0=0$, $i=1,...,n$;
\item[({\bf ii})] \begin{itemize}\item[$\bullet $] $\e_1=0$, 
                          \item[$\bullet $] $\e_2\in \spn e_1,e_2\es$, 
                          \item[$\bullet $] $\e_3\in \spn e_1,e_2,e_3\es;$
           \end{itemize}
\item[({\bf iii})] $\sum_{i=1}^n \e_i=0.$
\end{itemize}
	Conditions ({\bf ii}) above may be regarded as infinitesimal ``slice"
conditions while ({\bf iii}) is the infinitesimal closing condition.  Let us
write each of the component vectors of $\e$ as a sum of vectors
$\dl_i+\mu_i=\e_i$, where $\mu_i$ is the projection of $\e_i$ onto $\spn e_4,e_5\es$ for each $i$.  
Then, conditions ({\bf ii}) above imply that any $\e$ in $T_{U_0}S_0$ has the form
$$\e=(0,\dl_2,\dl_3,\dl_4+\mu_4,...,\dl_n+\mu_n).$$
At the same time, condition ({\bf i}) implies $\dl_2$ has $1$ degree of freedom 
within $\spn e_1,e_2\es$ whereas the remaining $\dl_i$, $i=3,...,n$, have
$2$ degrees of freedom within $\spn e_1,e_2,e_3\es$.  Finally, condition ({\bf iii}) says that
\begin{eqnarray}
\sum_{i=2}^n\dl_i=0 &\Leftrightarrow& -(\dl_2+...+\dl_{n-1})=\dl_n\in \spn e_1,e_2,e_3\es, \nonumber \\
\sum_{i=3}^n\mu_i=0 &\Leftrightarrow& -(\mu_3+...+\mu_{n-1})=\mu_n\in \spn e_4,e_5\es. \nonumber
\end{eqnarray}
	Now, the linear action of $H_2$ is clearly trivial on the $(2n-6)$-dimensional subspace of $T_{U_0}S_0$ spanned by the $\dl_i, \ i=2,...,n-1$.  However, it is a standard diagonal circle action on each of the $n-4$, two-dimensional subspaces spanned by each pair of linearly independent $\mu_i, \ i=1,...,n-1$.  The proposition follows. $\sq$
	
	The local structure of $\Mr$ at a singularity of type $3$ can be unveiled with a similar argument, as outlined below.  Let $Q_0$ be a degenerate polygon of type $3$.  As for the case of type 2 singularities, we may, if necessary, relabel the edges of $Q_0$ so that the first two span a 2-dimensional subspace of $\R^5$.  We may also choose a canonical representative $V_0=(e_1,v_2^0...,v_n^0)\in (S^4)^n$ of $[Q_0]$, where $v_2^0\in \spn e_1,e_2\es^+$, and define a $(4n-12)$-dimensional slice through $V_0$ for the $SO(5)$-action on $\Ur$ by
\begin{eqnarray}
S_0=\{V=(v_1,...,v_n)\in (S^4)^n: v_1&=&e_1, v_2\in\spn e_1,e_2\es, \nonumber \\
                             r\cdot V&=&\sum_{i=1}^n r_iv_i=0\}. \nonumber
\end{eqnarray}
	The smooth nature of $S_0$ near $V_0$ allows us to linearize our problem by considering the induced action of $H_3\simeq SO(3)$ on $T_{V_0}S_0$ given by the isotropy representation.  An argument parallel to the one offered in the previous proposition, shows that an arbitrary element $\e$ of $T_{V_0}S_0$ may be written as $\e=(0,\g_2,\g_3+\eta_3,...,\g_n+\eta_n)$, where $\g_i$, $i=2,...,n$, is the projection of $\e_i$ onto $\spn e_1,e_2\es$ and has 1-degree of freedom within this span.  Similarly, $\eta_i, \ i=3,...,n$, has $3$-degrees of freedom within $\spn e_3,e_4,e_5\es$. Taking into account the infinitesimal closing condition, we see that the fixed point set of $H_3$ is equivalent to $\R^{n-3}$.  We obtain the following proposition.

\begin{prop}
If $[Q_0]$ is a type $3$ singularity, then there exists a neighborhood of $[Q_0]$ in $\Mr$ isomorphic to a neighborhood of $0$ in 
$$\R^{n-3} \times [(\R^3)^{n-3}/SO(3)].$$
The linear factor $\R^{n-3}$ corresponds to planar deformations of $[Q_0]$.  $\sq$
\end{prop}

	We finish this section with a geometric description of a singularity of type $2$ for the moduli space of $5$-gons and $6$-gons, respectively.  For $n=5$, Proposition $2.6$ tells us that any type $2$ singularity has a neighborhood in $\Mr$ which is analytically isomorphic to a neighborhood of $0$ in $\R^4\times \R^{>0}$.  That is, $\Mr$ looks like a $5$-dimensional smooth manifold with boundary near this type of singularity.  The $4$-dimensional boundary component corresponds to the smooth, ``horizontal" deformations, namely those generated by the linearly-independent vector fields which locally span the singular locus $D_r$.  (Note that indeed, $D_r\simeq \R^4$ near any type 2 singularity as, in this case, no infinitesimal deformations lead to planar polygons).  The fifth ``transversal" direction corresponds, of course, to the non-degenerate deformations	of the polygon.	

	For $n=6$, any type $2$ singularity possesses a neighborhood in $\Mr$ isomorphic to $\R^6\times [(\R^2\times \R^2)/SO(2)]$.  The first component corresponds to the infinitesimal deformations along the (locally) smooth 6-dimensional singular locus $D_r$.  The transversal component is, in this case, equivalent to a 3-dimensional homogeneous quadratic cone.  Indeed, the linear action of $H_2\simeq SO(2)$ along transverse directions,
$$
H_2\times(\R^2\times \R^2)\rightarrow(\R^2\times \R^2):(h,(X,Y))\mapsto (hX,h^{-1}Y),
$$
induces an action on the polynomial ring $\R[X,Y]$, $X=(x_1,x_2)$, $Y=(y_1,y_2)$ given by
$$
H_2\times \R[X,Y]\rightarrow\R[X,Y]:(h,f(X,Y))\mapsto f(hX,h^{-1}Y).
$$
It is then immediate that $p_1=x_1^2+x_2^2$, $p_2=y_1^2+y_2^2$, $p_3=x_1y_1-x_2y_2$ and $p_4=x_2y_1+x_1y_2$ lie in the ring of invariant functions, $\R[X,Y]^{H_2}$, and satisfy the relation $p_1p_2=p_3^2+p_4^2$.  Let $V_p$ denote the (irreducible) $3$-dimensional semi-algebraic variety associated to the ring $\R[p_1,p_2,p_3,p_4]/\langle{p_1p_2-p_3^2+p_4^2}\rangle$, $p_1, p_2\geq 0$.  Then $V_p$ is isomorphic to $\R[X,Y]^{H_2}$.  Indeed, any given $(p_1,p_2,p_3,p_4)$ in $V_p$ with $p_1\neq 0$, is the image of the $H_2$-orbit of $(X,Y)=(\sqrt{p_1},0,\frac{p_3}{\sqrt{p_1}},\frac{p_4}{\sqrt{p_1}})$, so that the standard map $\R[X,Y]^{H_2}\rightarrow V_p$ is onto.  Also, if two $H_2$-orbits map to a single $(p_1,p_2,p_3,p_4)$ in $V_p$, then the specific form of the invariant functions guarantees the existence of $h\in H_2$ mapping one orbit to the other, that is, the orbits are the same and the map is one-to-one.  It follows that, in the transversal direction, $\Mr$ is equivalent to a homogeneous quadratic cone cut out by the equation $p_1p_2=p_3^2+p_4^2$.

	Similar arguments can be applied to reveal the geometric characteristics of a type 3 singularity.  For instance, for $n=5$ the $H_3\simeq SO(3)$-action along transversal directions gives once more the usual $H_3$-action on $\R[X,Y]$, $X=(x_1,x_2,x_3)$, $Y=(y_1,y_2,y_3)$.  It then follows that $|X|^2$, $|Y|^2$, $X\cdot Y$, $|X\times Y|$ are invariant functions subject to the relation $(X\cdot Y)^2+|X\times Y|^2=|X|^2|Y|^3$.	
\bigskip

\section{Real analytic equivalence between $\Mr$ and the weighted quotient of $(\HP^1)^n$ by $PSL(2,\H)$}

In this section, we construct the weighted quotient of the configuration space of $n$ points in $\HP^1$ by $\Psl$, denoted $Q_{st}$, and exhibit a real analytic equivalence between $\Mr$ and $Q_{st}$.  Establishing this equivalence involves extending the natural action of $\Psl$ on the rank-one symmetric space $\Psl/SO(5)$ to its geometric boundary in a manner that is consistent with the $SO(5)$-action on $S^4$.  In this sense, this construction can perhaps be seen as a special case of the results found in \cite{LeebM}.

	We begin by recalling that the real Lie group $GL(2,\H)$ may be viewed as the space of $2\times 2$ invertible matrices with entries in the skew field of quaternions, $\H$.  The invertibility condition is equivalent to the requirement that the Dieudonn\'e determinant, $D(A)$, of any $A$ in $GL(2,\H)$ be non-zero.  If we impose the additional condition that $D(A)=1$, we obtain the (semi-simple) Lie group $SL(2,\H)$.  Then, $\Psl=SL(2,\H)/\pm 1$.

	Let $\HP^1\simeq S^4$, denote the quotient space $\H^2\setminus\{(0,0)\}/\H^*$, where we assume the group of units $\H^*$ acts by right multiplication on $\H^2$.  We define a left action of $\Psl$ on $\HP^1$ by linear fractional transformations as follows:
$$\Psl\times \HP^1\rightarrow\HP^1:\ g\cdot[q_1:q_2]\mapsto[aq_1+bq_2:cq_1+dq_2],$$
$$\mbox{where }g=\left( \begin{array}{cc}
                          a & b\\
                          c & d
                         \end{array}\right)\, \in \Psl \mbox{, }(q_1,q_2)\in\H^2\setminus\{(0,0)\}.$$ 

	Let $M\subseteq(\HP^1)^n$ be the collection of $n$-tuples of distinct points.  Then the left diagonal action of $\Psl$ on $(\HP^1)^n$ induces an action on $M$ such that the quotient space $M/\Psl$ is a Hausdorff, real manifold.

	Consider an $n$-tuple of positive real numbers, $r=(r_1,...,r_2)$.  As $\Mr\simeq \M_{\l r}$ for all  $\l\in \R^+$, we may adopt the normalization $\sum_{i=1}^n{r_i}=2$.
\begin{df}
A point $p=(p_1,...,p_n)\in(\HP^1)^n$ is called stable if 
$$\sum_{p_i=q}r_i<1$$
for all $q\in \HP^1$.
\end{df}

	It is immediate that if we restrict ourselves to (admissible) $r$ satisfying the non-degeneracy condition (see Section 2), all $p\in(\HP^1)^n$ are stable. In virtue of the results of Section 2, we shall adopt this assumption for the remainder of this section.  However, it is not hard to define semi-stable and nice semi-stable configurations as in \cite{LeebM}.

	Let $M_{st}$ denote the space of all stable points in $(\HP^1)^n$ and set $Q_{st}=M_{st}/\Psl$ with the quotient topology.  Recall from Section 2 that $\Ur\subseteq (S^4)^n$ represents the collection of closed $n$-gons, $\Pn$.  The first step in establishing the correspondence between $Q_{st}$ and $\Mr=\Ur/SO(5)$ is to note that the closing condition implies the stability condition in Definition 3.1.  As $SO(5)$ is a maximal compact subgroup of $\Psl$, the inclusion $\Ur\subseteq\M_{st}$ induces the (injective) quotient map
$$\xi:\Mr=\Ur/SO(5)\rightarrow Q_{st}=M_{st}/\Psl.$$

\begin{Th}
The quotient map $\xi$ gives a real analytic isomorphism between $\Mr$ and $Q_{st}$. 
\end{Th}

	Proving Theorem 3.2 essentially amounts to proving Lemma 3.4 below, which, modulo some preliminary results on the action of $\Psl$ on a space of discrete stable measures on $S^4$, guarantees the surjectivity of the map $\xi$.  The same ideas were used by Kapovich and Millson in \cite{KM} to establish a complex analytic equivalence between $\3$ and the weighted quotient of $(S^2)^n$  by $PSL(2,\C)$, constructed by Deligne and Mostow in \cite{DM}.  As it is also the case in this section, the arguments used in \cite{KM} to establish such equivalence rely on some special properties of the conformal barycenter constructed by Douady and Earle in \cite{DE} for stable measures on the sphere.

\begin{df}
A probability measure on $S^4$ is said to be stable if the mass of any atom is strictly less than $\frac{1}{2}$.
\end{df}

	For our choice of $r$, each point in $\Mr$ gives rise to a stable, finite probability measure of total mass 1 defined by
$$\mu=\frac{1}{2}\sum_{i=1}^nr_i\dl_{u_i},$$
where $\dl_{u_i}$ denotes the delta function on $S^4$ centered at $u_i$.
The center of mass, $C(\mu)$ of such a measure is given by
$$C(\mu)=\frac{\sum_{i=1}^nr_iu_i}{\sum_{i=1}^nr_i}=\frac{1}{2}\sum_{i=1}^nr_iu_i.$$
Note that the $\Psl$-action on $S^4$ by linear fractional transformations induces an action on probability measures on $S^4$ defined by
\begin{equation}
g\cdot\mu(A)=g_*\mu(A):=\mu(g^{-1}(A)),\mbox{ where } A\subseteq S^4 \mbox{ is a Borel set}.
\label{e:push}
\end{equation}

\begin{lem}
For each finite stable measure $\mu$ on $S^4$ there exists $g\in\Psl$ such 
that the center of mass $C(g\cdot \mu)=0$.
\end{lem}

\pf Consider the Cartan decomposition of $\Psl=SO(5) P,$
$$P=\left\{ \left( \begin{array}{cc}
                          \rho_1      &  q     \\
                          \bar{q} & \rho_2
              \end{array}\right),\ \rho_1\rho_2=1+|q|^2,\ \rho_1>0,\ \rho_2>0,\ q\in\H\right\},$$
and the natural action of $\Psl$ on $P$ given by
$$\Psl\times P\rightarrow P:(g, x)\mapsto g x g^*, \ g^*=\bar{g}^t.$$
We will first show that the geometric boundary of $P$ may be identified with $S^4$ and that the above $\Psl$-action may be continuously extended to an action on this geometric boundary which  coincides with the action on $S^4$ by linear fractional transformations.
To this purpose, we identify the vector space $P$ with $\Rup=\{(x_1,x_2,x_3,x_4,x_5), \ x_5>0\}$ by setting
\begin{eqnarray}
x_i&=&\frac{q_i}{\rho_2},\mbox{ }i=1,2,3,4; \nonumber \\
x_5&=&\frac{1}{\rho_2}>0 \nonumber
\end{eqnarray}
where $q=q_1+iq_2+jq_3+kq_4\in \H$.  The $\Psl$-action can be written in terms of the $x_i$ coordinates as
\begin{eqnarray}
v_x & \mapsto & \frac{{|x|^2}a\bar{c}+b\bar{v_x}\bar{c}+av_x\bar{d}+b\bar{d}}{{|x|^2}|c|^2+d\bar{v_x}\bar{c}+cv_x\bar{d}+|d|^2}\nonumber \\
x_5 & \mapsto &  \frac{x_5}{{|x|^2}|c|^2+d\bar{v_x}\bar{c}+cv_x\bar{d}+|d|^2}>0\nonumber
\end{eqnarray}
where $|x|^2=\sum_{i=1}^5{x_i^2}$, \ $v_x=x_1+ix_2+jx_3+kx_4\in \H, \ g=\left( \begin{array}{cc}
                                          a & b\\
                                          c & d
                                         \end{array}\right)$.  Let us now continuously extend this $\Psl$-action to the geometric boundary of $\Rup$ consisting of $\{x\in \R^5: x_5=0\}\cup\{\infty\}$ by 
\begin{eqnarray}
g\cdot (x_1,x_2,x_3,x_4,0)&:=&\lim_{x_5\rightarrow 0} g\cdot  (x_1,x_2,x_3,x_4,x_5)\nonumber \\
                          &= &\lim_{x_5\rightarrow 0} \frac{{|x|^2}a\bar{c}+b\bar{v_x}\bar{c}+av_x\bar{d}+b\bar{d}}{{|x|^2}|c|^2+d\bar{v_x}\bar{c}+cv_x\bar{d}+|d|^2}\nonumber \\
                          &= &\frac{|v_x|^2a\bar{c}+b\bar{v_x} \bar{c}+av_x\bar{d}+b\bar{d}}{|v_x|^2|c|^2+d\bar{v_x}\bar{c}+cv_x\bar{d}+|d|^2}\nonumber \\
                          &= &(av_x+b)(cv_x+d)^{-1}\nonumber.
\end{eqnarray}
Upon identifying each $v_x\in\H$ with the corresponding $[v_x:1]\in\HP^1$, we obtain an action of $\Psl$ on $\H\cup{\infty}=\HP^1$ which coincides with the original $\Psl$-action on $\HP^1$ by linear fractional transformations.
	
	Let $B^5$ denote the closed unit ball in $\R^5$.  One can introduce coordinates $y_i,\ i=1,...,5$ in $\R^5$, and define an invertible map $\R^5_{> 0}\rightarrow B^5$ by
\begin{eqnarray}
y_i & = & \frac{2x_i}{1+|x|^2},\ i=1,...,4;\nonumber \\
y_5 & = & \frac{1-|x|^2}{1+|x|^2}. \nonumber
\end{eqnarray}
If we then express the $\Psl$-action in terms of the $y$-coordinates for the ball model, we may directly check that the map sending a finite stable measure $\mu$ on $S^4$ to its center of mass $C(\mu)\in B^5$ fails to be $\Psl$-equivariant.  However, in their paper \cite{DE}, Douady and Earle construct a conformal barycenter $B\in B^n$ associated to every stable probability measure on $S^{n-1}$.  For $n=5$, the conformality of $B$ amounts to the $\Psl$-equivariance of the assignment $\mu\mapsto B(\mu)$, namely,
$$B(g\cdot\mu)=g(B(\mu)),$$
where $\Psl$ acts on $B^5$ as indicated above, and on measures by push-forward, as defined in (\ref{e:push}).  It turns out that, for stable $\mu$, the conformal barycenter $B(\mu)$ coincides with the unique zero of a vector field $\xi_{\mu}$ defined on (the appropriate) $B^n$.  For $n=5$ and stable, finite $\mu$, we have that
$$\xi_\mu(y)=\frac{1}{2}\sum_{i=1}^n\left({\frac{1-|y|^2}{|y-u_i|^2}}\right)^4r_i (y-u_i),$$
where $y\in B^5$, $u_i\in S^4$ for all $i$, and $\mu$ is defined by $u=(u_1,...,u_n)$.  As $\xi_\mu(0)=\frac{1}{2}\sum_{i=1}^n r_iu_i$, one immediately verifies that: 
$$B(\mu)=0\Leftrightarrow C(\mu)=0.$$
Lemma 3.4 then follows at once from the transitivity of the $\Psl$-action on $B^5$. $\sq$  

	The real analytic equivalence in Theorem 3.2 is thus established. $\sq$

\bigskip

\section{Gel{\cprime}fand-MacPherson correspondence over the quaternions}

The classical Gel{\cprime}fand-MacPherson correspondence \cite{GMacP} asserts that the quotient of the generic part of $Gr_{\R}(p, q)$ by the Cartan subgroup of $PGL(p+q, \R)$ is diffeomorphic to the space of equivalence classes of generic configurations of $(p+q)$ points in $\RP^{q-1}$. Later on, it was realized (see e.g. \cite{Kapr}) that in the complex context, there is an isomorphism between appropriately chosen symplectic quotients of these spaces, as well as between the corresponding GIT quotients. 

	The goal of this section is to show that a similar result holds over the quaternions. More precisely, let $n=p+q$, and let $P$ be the subgroup of $SL(n, \H)$ which preserves a $q$-dimensional subspace of $\H^n$. Then $Gr(p,q)=SL(n, \H)/P\simeq Sp(n)/Sp(p)Sp(q)$ is the quaternionic Grassmannian with the $PSL(n,\H)$-action.  In this case, the maximal compact subgroup of $SL(n, \H)$ is $Sp(n)$, and $\S^n\simeq Sp(1)^n\subset Sp(n)$ is the diagonal subgroup of $Sp(n)$ and thus, a natural analog of the torus.  We can then consider the tri-momentum map $\mu: Gr(p,q)\to \R^n$ defined similarly to the moment map in the complex case \cite{F1}.  The level sets of $\mu$ are $\S^n$-invariant, so for any $x=(x_1, ..., x_n)\in \R^n$, one gets the reduced space $X_x:=\mu^{-1}(x)/\S^n$, which is not smooth in general, as the action of $\S^n$ on $\mu^{-1}(x)$ is free only generically. 

	On the other hand, let us consider the space ${\tilde Y}=(\HP^{p-1})^n$ with the diagonal $PSL(p, \H)$-action. There is a map $\phi$ from ${\tilde Y}$ to the (real) projectivization of $\calH_p$ ---the space of quaternionic Hermitean $p\times p$ matrices, denoted by $\P\calH_p$. In terms of the homogeneous coordinates $[w_1^{(i)}:\cdots : w_p^{(i)}]$ on the $i^{\mbox{{\small th}}}$ component $\HP^{p-1}$, the map $\phi$ sends a point with these coordinates to the projectivization of the matrix whose $(i,j)^{\mbox{{\small th}}}$ entry is given by $\sum_{l=1}^p w_l^{(i)}{\bar w}_l^{(j)}$. There is an obvious $Sp(p)$-action on $\P\calH_p$ by conjugation, and if $\calO$ is an orbit of this action, then $\phi^{-1}(\calO)$ is preserved by $Sp(p)$. We have the following:
\begin{prop}
For the appropriate choices of $x$ and ${\cal O}$, the corresponding quotients $X_x$ and $Y_{\calO}:=\phi^{-1}(\calO)/Sp(p)$ are homeomorphic, while their open dense smooth parts are in fact diffeomorphic. 
\end{prop}

\pf
For both quotients, we can start with the same space $\H^{pn}$, with coordinates $(w_i^{(j)})$, $1\le i\le p$, $1\le j\le n$, from where we have two maps $\mu: \H^{pn}\to \R^n$ given by $x_i=\sum_{j=1}^n |w_i^{(j)}|^2$ and $\phi: \H^{pn}\to \calH_p$ given as above, prior to the projectivization.  One can easily establish that the level sets of the first map are preserved by the action of $Sp(p)$ and the level sets of $\phi$ are preserved by $\S^n$.  Our result then follows from a simple application of reduction in stages and the fact that $Gr(p,q)$ is obtained from $\H^{pq}$ by choosing a level set of $\phi$ and quotiening by $Sp(p)$, while $(\HP^{p-1})^n$ is obtained by choosing a level set of $\mu$ and quotiening by $\S^n$. Therefore, we see that for the appropriate choices of $x$ and ${\cal O}$ we obtain a bijection on the level of orbits.  $\sq$
\bigskip

An alternate way of establishing the quaternionic version of the
Gel{\cprime}fand-MacPherson correspondence by using a certain natural
involution will be presented in Section 7. 
          
\section{Gel{\cprime}fand-Tsetlin coordinates and application to $\M_r$}
 
The relationship between the bending flows on the moduli spaces of polygons in $\R^3$ defined by Kapovich and Millson in \cite{KM} and the Gel{\cprime}fand-Tsetlin (GT) system on the Grassmannian $Gr_{\C}(2, n)$ was exhibited by Hausmann and Knutson in \cite{HK}. In this section we simply follow the road paved in \cite{HK} to show that one can obtain a similar relation in the quaternionic
context. 

	First of all, let us recall that for any $n\times n$ quaternionic Hermitean matrix $A\in \calH_n$, the set of $n$ real eigenvalues (and hence all real invariant polynomials) is well-defined by formulas analogous to the complex case. (Note that this is not the case for a general matrix with quaternionic entries.) Therefore, the whole system of Gel{\cprime}fand-Tsetlin coordinates $\{ \l_i^{(j)}\}$ described below, is also well-defined on $\calH_n$.  Recall that given a matrix $A\in \calH_n$, we may denote by $A^{(j)}$ the $j\times j^{\mbox{{\small th}}}$ upper-left corner submatrix of $A$. Then, by definition, $\l_i^{(j)}(A)$ is the $i^{\mbox{{\small th}}}$ eigenvalue of $A^{(j)}$ (these are arranged in non-increasing order). The whole system $\{ \l_i^{(j)}\}$ satisfies the interlacing property: 
$$\l_i^{(j)}\ge \l_i^{(j-1)}\ge \l_{i+1}^{(j)}.$$

	Due to the collision of eigenvalues, the functions $\l_i^{(j)}$ are not smooth on all of $\calH_n$, but only continuous.  However, their importance in the complex case, as was shown by Guillemin and Sternberg \cite{GSGT} \cite{GST} based on Thimm \cite{Th}, is manifested by the fact that, generically, these functions form a set of complete integrals on the co-adjoint orbits of $U(n)$ (that is, on the level sets of the functions $\l_i^{(n)}$) for the Lie Poisson structure. Moreover, generically, the level sets of the rest of the (GT) coordinates are Lagrangian tori. 

	Similarly, one can define (GT) coordinates on $\calH_n$ in the quaternionic case, as every element of $\calH_m$ has $m$ real eigenvalues.  In this case, the orbits of the natural $Sp(n)$-action by conjugation are analogues of the co-adjoint orbits. To see that the interlacing property holds, one can use the embedding $\nu: \calH_n\hookrightarrow \calH_{2n}^{\C}$ defined in Equation (\ref{e:nu}) in the next section, and use the interlacing property for the eigenvalues of the complex Hermitean sub-matrices. The main point is that each eigenvalue in the embedding defined by $\nu$ will appear twice.   

	One can also see that the generic level sets of these functions on the orbits $\calO$ of the $Sp(n)$-action are diffeomorphic to $\S^m$, the diagonal subgroup of $Sp(n)$, for $m=\dim (\calO)/4$.  This fact is also easy to prove by an inductive argument.

	Let $v_i, \ i=1,...,n$, denote the vertices of a arbitrary $n$-gon in $\Mr$ and $d_i=v_{i+2}-v_1, \ i=1,..., n-3$, its $n-3$ diagonals.  Set $\ell_i=||d_i||,\ i=1,...,n-3$.  We would like to show that an analogue of Theorem 5.2 in \cite{HK} holds in our case. Namely, that if one fixes a level set of the diagonal lengths, $L_\ell$, for a generic polygon in $\Mr$, then the transitive spheroid action on $L_{\bf \ell}$ comes from the residual action of $\S^{2n-4}$ on a level set of the (GT) system on $Gr(2, n)$.  Indeed, an element of $Gr(2,n)$ can be represented by a $2\times n$ matrix $M$ with quaternionic entries, where the columns span a $2$-dimensional subspace of $\H^n$. If $(a_1,..., a_n)^t$ and $(b_1, ..., b_n)^t$ are the columns, then we can also introduce the truncated $2\times i$ submatrices $M_i$, where only the first $i$ rows remain.  As in Section 3, let $M^*=\overline{M}^t$.  Then the non-zero eigenvalues of the $2\times 2$ quaternionic Hermitean matrix 
$$
M_i^*M_i=\sum_{j=1}^i
\left(
\begin{array}{cc}
|a_j|^2 & {\bar a}_jb_j \\
{\bar b}_j a_j & |b_j|^2
\end{array} \right)
$$ 
are the same as the non-zero eigenvalues of the matrix $M_iM_i^*\in\calH_i$.  The proof of this fact is analogous to Section 5 of \cite{HK} and gives the result. 
\bigskip

\section{Relations with complex geometry}

Let us identify $\C^{2m}$ with $\H^m$ as follows. The point $(z_1, ..., z_{2m})\in \C^{2m}$ 
corresponds to the point $(q_1, ..., q_m)\in\H^m$ if $q_i=z_{2i-1}+{\bf j}z_{2i}$ for 
$1\le i\le m$. Using this identification, let $J$ be the real operator on $\C^{2m}$ which comes from the right
multiplication by ${\bf j}$ on the space $\H^m$. Since $J^2=-\mbox{Id}$, the action 
of $J$ on $\C^{2m}$ extends to an involution $\theta$ on all complex partial flag manifolds. 
However, this involution is only real and not a complex diffeomorphism. If all the dimensions 
of the subspaces are even, then the fixed point set of $\theta$ is clearly the quaternionic 
(partial) flag manifold: 
$$
(F_{2m}^{\C}(2m_1, ..., 2m_k))^{\theta}=F_m^{\H}(m_1, ..., m_k).
$$
	Moreover, if $\omega$ is an invariant K\"ahler form on the complex flag 
manifold, then the fixed point set of $\theta$ is a Lagrangian submanifold with 
respect to $\omega$. Let us now restrict 
our attention to the case of the complex Grassmannian $Y=Gr_{\C}(2,4)$ of complex 
2-planes in a 4-dimensional complex space.  In particular, the above discussion 
implies that $Y^\theta=\HP^1$.  We can always view the group $SL(n, \H)$ as a real 
form of the complex semi-simple group $SL(2n, \C)$, with the corresponding Satake 
diagram \cite{Araki} having odd numbered vertices painted and no arrows.  
The multiplicity of each restricted root is $4$.  In terms of matrices, we have the embedding
\begin{equation}
\nu: SL(n, \H)\hookrightarrow SL(2n, \C), \ 
A+B{\bf j}\ \mapsto \left(
\begin{array}{cc}
A & B \\ -{\bar B} & {\bar A}
\end{array}
\right).
\label{e:nu}
\end{equation}
If we let $J$ be the $2n\times 2n$ matrix
$$
\left(
\begin{array}{cc}
0 & 1 \\ -1 & 0
\end{array}
\right),
$$
then one can define the involution $\t$ on $SL(2n, \C)$ by $\t(C)=-J{\bar C}J$, which defines the real form $SL(n, \H)$. One can see that using the notation $\t$ for both the involution just defined and the involution on $\C^{2n}$ presents no harm, since both have the same origin and the action of $SL(n, \H)$ on the fixed point set of $\t$ in flag manifolds comes from the action of $SL(2n, \C)$.  

	Let us now consider the $n$-fold product $Y^n=\underbrace{Y\times\cdots\times Y}_n$ with the diagonal action of the group $G=SL(4, \C)$. We will always assume that $n>4$, since we would like generic orbits to depend on continuous parameters.  The space $Y^n$ has two different interpretations. On one hand, it can be viewed as the configuration space of $n$-tuples of ordered points on $Y=Gr_{\C}(2,4)$.  Another point of view is that the space $Y^n$ is the configuration space of ordered lines in $\CP^3$.  In both situations the group $G$ acts by identifying projectively equivalent configurations.  Consider a weighted configuration of lines $\{ l_1, ..., l_n\}$ and corresponding weights $(r_1, ..., r_n)$ subject to the normalization $\sum_{i=1}^n r_i=2$.  Then, following \cite{GIT}, pages 86-88, we say that such an $n$-tuple of 
lines in $\CP^3$ is (properly) {\it stable} if the following three conditions hold:
\begin{itemize}
\item[({\bf i})] For every point $p\in \CP^3$, the sum of weights of lines passing through 
$p$ is $< 1$.
\item[({\bf ii})] For every line $l\in\CP^3$, the sum of weights of lines intersecting 
$l$ plus twice the sum of weights of lines coincident with $l$ is $< 2$.
\item[({\bf iii})] For every plane $\Pi\subset\CP^3$, the sum of weights of lines 
contained in $\Pi$ is $< 1$.
\end{itemize}
As usual, for a semi-stable configuration, one replaces $<$ by $\le$ in the
above definition.  
In particular, if a $G$-orbit contains a stable point, then it is entirely stable.  
Further, stable orbits are $15$-dimensional, and hence, as large as possible. 
We notice that there always exist semi-stable configurations, which are not stable. 
The easiest way to see this is to choose a semi-stable configuration of lines each 
of which meets a certain fixed line in $\C\P^3$ at a point. Then, in ({\bf ii}) 
above, we will have a strict equality.  

	We would now like to consider the involution $\t$ on $Y^n$ and show that it 
descends to GIT quotients by the action of $G$.  We will see that the fixed point 
sets of this action on the quotients are related to our polygon spaces $\M_r$.  
More precisely, on the Grassmannian $Y$ there is an essentially unique line bundle 
$\L$, dual to the second exterior power of the tautological plane bundle.  One can 
also view $\L$ as the pull-back of the canonical line bundle ${\cal O}(1)$ over 
projective space, in the Pl\"ucker embedding of $Y$.  Therefore, all the choices 
for linearization on the product $Y^n$ are given by taking the tensor product of 
the pull-backs of the line bundles $\L^{\otimes a_1}, ..., \L^{\otimes a_n}$.  
We assume that all $a_i >0$ so that the corresponding line bundle is ample.  
Let $a=(a_1, ..., a_n)$ and let us set $r_i=2a_i/\sum_{i=1}^n a_i$.  Whenever we 
need to specify the linearizing line bundle, we will denote the corresponding GIT 
quotient by $(Y^n/G)_r$, with $r=(r_1, ..., r_n)$.  Now we would like to show that 
the action of $\t$ on $Y^n$ maps orbits of $G$ to orbits of $G$. Indeed, if we let 
$A$ be an element of $G$, a simple computation shows that $\t$ maps the element 
$A\cdot x$ to the element $\t(A)\cdot\t(x)$. In particular, if $x$ is in the fixed 
point set of $\t$ and $A\in SL(2, \H)$, then $\t$ will stabilize $A\cdot x$. We also
note that $\theta$ maps semistable configurations to semistable configurations. 
This immediately shows that the action of $\t$ descends to the GIT quotients 
$(Y^n/G)_r$. 

	Next, we give an interpretation of the quaternionic Gel{\cprime}fand-MacPherson 
correspondence using the involution $\t$.  In complex algebraic geometry, for an 
$n$-tuple of positive integers $a=(a_1, ..., a_n)$, we have the following isomorphism 
of GIT quotients:
\begin{equation}
Gr_{\C}(2,4)^n/SL(4,\C)\simeq Gr_{\C}(4,2n)/GL(2, \C)^n,
\label{e:GM2}
\end{equation}
where the corresponding linearizing line bundles are defined as follows.  On the left hand side, the linearizing line bundle is the same as before, the tensor product of the pull-backs of $\L^{\otimes a_i}$. On the right hand side, we think of $GL(2, \C)^n$ as the block-diagonal subgroup of $GL(2n, \C)$. Notice that the one-dimensional center of the latter group is a subgroup of the former and acts trivially on $Gr_{\C}(4, 2n)$.  In particular, the dimension count is correct.  The linearization of the action of $GL(2, \C)^n$ on the right can be defined analogously to (2.4.3) of \cite{Kapr}. More precisely, let $Z_i$ be the center of the $i^{\mbox{{\small th}}}$ copy of $GL(2, \C)$ in the $n$-fold product.  The subgroup $Z_i$ is identified with $\C^*$ using our matrix embedding.  Our point is that both varieties in Equation (\ref{e:GM2}) are the same as the projective spectrum of the ring 
\begin{equation}
R=\oplus R_d,
\label{e:R_d}
\end{equation}
where $R_d$ consists of polynomials $\Phi(M)$ in entries of a $(4\times 2n)$-matrix, such that
\begin{itemize}
\item[{\bf 1.}] $\Phi(gM)=\Phi(M)$ for the left action of $g\in SL(4,\C)$, 
\item[{\bf 2.}] $\Phi(M\cdot h)=\chi(h)\Phi(M)=t^{da}\Phi(M)$, where $h\in GL(2,\C)^n$.
\end{itemize}
Here, $\chi$ is the character of $GL(2,\C)^n$ corresponding to the character $t^{da}$ of its center, where $t\in Z_1\times\cdots\times Z_n\simeq (\C^*)^n$, and $t^a=t_1^{a_1}\cdots t_n^{a_n}$ component-wise.  The correspondence above can be easily extended to an equivalence between GIT quotients of Grassmannians $Gr_{\C}(k,d)^n/SL(d, \C)$ and $Gr_{\C}(d, nk)/GL(k, \C)^n$, and even to an equivalence between more general partial flag manifolds (for appropriate choices of linearizations).  Now, one can see that the correspondence defined by Equation (\ref{e:GM2}) is preserved by the action of the involution $\t$.  This yields an alternate way of establishing the Gel{\cprime}fand-MacPherson correspondence in the quaternionic context.  We summarize our discussion in the following proposition.

\begin{prop} For an admissible $n$-tuple of weights $r=(r_1, ..., r_n)$ 
the GIT quotient $(Y^n/G)_r$ is a projective variety of dimension $4n-15$, 
determined by the projective spectrum of the ring $R$, defined in 
Equation (\ref{e:R_d}) above. The GIT quotients 
$$Gr_{\C}(2,4)^n/SL(4, \C)\ \ \ \ {\rm and} \ \ \ \ Gr_{\C}(4,2n)/GL(2, \C)^n$$ 
are isomorphic for appropriate choices of linearizations.  
The action of the involution $\t$ on $Y^n$ descends to the quotient $(Y^n/G)_r$. $\hfill\Box$
\end{prop}
 
	Let us now consider symplectic reductions of $Y^n$ with respect to the
$K=SU(4)$-action. The symplectic form on $Y^n$ is the sum of pull-backs of
$K$-invariant forms on each multiple.  These are proportional, with the
corresponding coefficients $r_i$, to the one determined by the line bundle
$\L$.  Let us identify $\su_4^*$ (the dual of the Lie algebra $\su_4$) with the
space of complex traceless Hermitean $4\times 4$ matrices.  Then, the fixed
point set of $\t$ is again given by the 
image of the space of traceless quaternionic $2\times 2$ Hermitean matrices
under the embedding defined by the map $\nu$ in Equation (\ref{e:nu}).  Since
the involution $\t$ maps unitary matrices to unitary matrices, then it also
maps $K$-orbits on $Y^n$ to $K$-orbits.  Therefore, $\t$ descends to the
symplectic quotient $(Y^n//K)_r$, which is homeomorphic to the GIT quotient
$(Y^n/G)_r$ as explained in \cite{GIT}, for example.  Our convention is that we
always reduce at the zero level set of the moment map, unless stated otherwise.
The quotient $(Y^n//K)_r$ is not expected to be smooth even for a sufficiently
general choice of weight vector $r$, as even then one can always produce a
configuration with a non-trivial connected stabilizer as explained above.  
However, we have the generic smoothness on the symplectic side, since the 
$K$-action is, generically, locally free. In any case, the quotient space 
$(Y^n//K)_r$ has the structure of a stratified symplectic manifold, in the 
terminology of \cite{SjamLer}. Thus, we arrive at the following result:

\begin{prop} For an admissible $n$-tuple of weights $r=(r_1, ..., r_n)$, 
the symplectic quotient $(Y^n//K)_r$ is a stratified symplectic manifold.  
The action of $\t$ on $Y^n$ descends to a smooth anti-symplectic involution 
on each smooth symplectic stratum of $(Y^n//K)_r$, where its fixed point set 
is a Lagrangian submanifold.
\end{prop}

\pf Recall that if $\tau$ is a smooth involution on a manifold $M$, then the 
fixed point set $M^{\tau}$ is a smooth submanifold of $M$. Indeed, one can use 
the fact that for each fixed point $p$ in $M^{\tau}$, the tangent space $T_pM$ 
splits as $V^+\oplus V^-$ according to the eigenvalues $1$ and $-1$ of the
induced 
tangent space map $\tau_*(p)$.  Then, one can use a $\tau$-invariant Riemannian 
metric on $M$ together with the exponential map on $T_pM$ to see that $p$ is a 
smooth point of $M^{\tau}$. The fact that $\t$ descends to the anti-symplectic 
involution on $(Y^n//K)_r$ immediately follows from the fact that $\t$ is an 
anti-symplectic involution on $Y^n$.  It follows from Proposition 2.3 of 
\cite{OSS} that the fixed point set of $\t$ on each smooth symplectic stratum,
is a Lagrangian submanifold. $\sq$
 
	Note that the dimension of the fixed point set of $\t$ on the big
smooth open subset of $(Y^n//K)_r$ is $4n-15$, half of the dimension of $(Y^n//K)_r$.  

	Let us now recall a general Lagrangian reduction procedure that was considered 
in Section 7 of \cite{OSS}.  One can interpret the space $(Y^n//K)_r$ as the space of 
solutions of the equation $A_1+\cdots +A_n=0$, where $A_i\in\su^*_4$ is in the orbit 
of ${\rm diag}(r_i, - r_i, r_i, -r_i)$, modulo the diagonal action of $SU(4)$.  
Our space of polygons $\M_r$ can also be interpreted as the space of solutions of the 
same equation $B_1+\cdots +B_n=0$ but with $B_i$ belonging to the space $\calH_2$ of 
quaternionic Hermitean traceless $2\times 2$ matrices, modulo the diagonal action of 
the group $Sp(2)$ by conjugation. To see that these two spaces are related by an 
involution, let us construct a map: 
$$
\psi: \M_r\to (Y^n//K)_r^{\t},
$$ 
and study its properties. The construction is already clear from the previous
paragraph.  Indeed, using the embedding $\nu$ defined in Equation (\ref{e:nu}),
one can map the $Sp(2)$-orbit of an $n$-tuple of matrices 
$(B_1$, ..., $B_n)\in\calH_2$ satisfying $B_1+\cdots+ B_n=0$ to the $SU(4)$-orbit 
of an $n$-tuple of matrices from $\su^*_4$ satisfying the same equation, where
we let $A_i=\nu(B_i)$.  The image of such orbit not only lies entirely in an 
$SU(4)$-orbit, but also belongs to the fixed point set of the involution $\t$.  
Therefore, when we pass to the quotients, we get a well-defined map $\psi$ defined as above. 

\begin{prop} The map $\psi$ is finite and generically four to one. 
\end{prop} 

\pf The finiteness of $\psi$ follows from \cite[Proposition 2.3 (iii)]{OSS},
which, applied to our case, implies that an $Sp(2)$-orbit through an $n$-tuple
$(B_1$,...,$B_n)\in\calH_2$ satisfying $B_1+\cdots +B_n=0$ is open in the intersection of  
the $SU(4)$-orbit through $(B_1, ..., B_n)$ with $\calH_2$. Assume that 
two different $Sp(2)$-orbits, through the points $x_1$ and $x_2$ from 
the zero level set are mapped to the same point by $\psi$. This would imply that
there exists an element $g\in SU(4)$ such that $\nu(x_1)=g.\nu(x_2)$. 
Now if we apply the involution $\theta$ defining the quaternionic subspace inside the 
complex one to both sides, we will see that $\nu(x_1)=\theta(\nu(x_1))$ 
$=\theta(g)\theta(\nu(x_2))=\theta(g)\nu(x_2)$, where $\theta$ applied to an element 
from $SU(4)$ stands for the involution inside $SU(4)$ defining $Sp(2)$. 
Since, generically, points have central stabilizer and
$\nu$ is an injection, we see that $g=c\cdot\theta(g)$, where $c$ is in the center of 
$SU(4)$, which consists of four elements. Thus, generically, each $SU(4)$ orbit 
will contain four $Sp(2)$-orbits, which are stabilized by $\theta$. $\sq$ 

	Notice that both $\M_r$ and $(Y^n//K)_r^{\t}$ are real analytic varieties 
of dimension $4n-15$. Thus the map $\psi$ surjects onto a connected component of 
$(Y^n//K)_r^{\t}$ cf. \cite[Corollary7.2]{OSS}, and the image of $\psi$ consists 
of those $SU(4)$-orbits which are preserved by $\t$ and intersect its fixed point set.  

However, the map $\psi$ is not surjective. The question of surjectivity, as shown 
in \cite{F3}, boils down to the question of how many conjugacy classes of involutions 
inner to $\theta$ there are in $PU(4)$. This fact holds in general 
for a compact group $K$ of adjoint type and an involution $\theta$. 
An involution $\tau$ is called {\it inner}
to $\theta$ if there exists an element $g$ from $PU(4)$ such that 
$\tau={\rm Ad}_{g}\circ\theta$. Two involutions $\tau$ and $\tau'$ are called 
conjugate, if there exists a group element $g$ such that 
$\tau'={\rm Ad}_{g}\circ\tau\circ{\rm Ad}_{g^{-1}}$. 
For the group $PU(4)$ and the involution $\theta$ defining the symlpectic subgroup, 
there is another conjugacy class of involutions represented by $\tau$, the complex 
conjugation, whose fixed point set is the group $PO(4)$. The corresponding involution,
also denoted by $\tau$ on $Y=Gr_{\C}(2,4)$ has the real grassmannian $Y^\tau=Gr_{\R}(2,4)$ 
as its fixed point set. Therefore, we can also consider the lagrangian quotient of 
the $n$-fold product of the real grassmanian of two-planes in $\R^4$ by the diagonal 
action of the group $PO(4)$. This quotient is defined as the quotient of 
the intersection of $Gr_{\R}(2,4)^n$ with the zero level set of the momentum map inside 
$Y^n=Gr_{\C}(2,4)^n$ and denoted by ${\mathcal T}_r$. 

Similar to the map $\psi$ defined above, we can define the map 
$$ \varphi: \ \ {\mathcal T}_r\to (Y^n//SU(4))^\tau_r.$$
It was proved in \cite{F3} that the spaces $(Y^n//SU(4))^\tau_r$ and 
$(Y^n//SU(4))^\theta_r$ are naturally identified, and the images of the maps 
$\psi$ and $\varphi$ are actually disjoint. It is easy to see that for an admissible
choice of $r$, the space ${\mathcal T}_r$ is non-empty (also of dimension $4n-15$).
The following follows from \cite{F3}: 

\begin{prop}
The space $(Y^n//SU(4))^\theta_r$ is 
homeomorphic to the disjoint union of $\psi({\mathcal M}_r)$ and $\phi({\mathcal T}_r)$. 
\end{prop}

	Now we would like to consider the Chow quotient $Y^n//G$, as defined by 
Kapranov in \cite{Kapr} and extend the action of $\t$ to an involution on $Y^n//G$, 
whose fixed point set $\M_n$  will have a similar meaning to the Grothendieck-Knudson 
space $\overline{M}_{0,n}$ ---a certain compactification of projective equivalence 
classes of $n$-tuples of distinct points on $\CP^1$.  In particular, we will construct 
surjective maps $\M_n\to\M_r$.  Thus, in a sense, the space $\M_n$ will serve as a 
``universal'' polygon space dominating all the $\M_r$ spaces at once.  First of all, 
let us fix a weight vector $r$ and consider the bi-rational morphism 
$\phi_r: Y^n//G\to(Y^n/G)_r$ defined in \cite{Kapr}. Now let us denote by $\M_{n,r}$ 
the fiber product $\M_r\times_{(Y^n/G)_r}Y^n//G$, where the map $\M_r\to (Y^n/G)_r$ 
is defined using the map $\psi$ from Proposition 6.3. Since the fixed point set of 
the action of $\t$ on the corresponding Chow variety is closed, the action of $\t$ 
extends to an involution on the Chow quotient $Y^n//G$.  There is natural map 
$\eta_r: \M_{n,r}\to (Y^n//G)^\t$, which is defined using the above construction of 
$\M_{n,r}$ and the fact that $\psi$ actually maps $\M_r$ to $(Y^n/G)_r^{\t}$. As in 
the above proposition, we can show that the map $\eta_r$ is generically injective.  
We define the space $\M_n$ as the common fiber product of the spaces $\M_{n,r}$ using 
the maps $\eta_r: \M_{n,r}\to(Y^n//G)^\t$. This product is actually finite for every 
$r$, because the isomorphism class of $\M_r$ is clearly determined by the combinatorial 
choice of $<0$, $>0$, or $=0$ for the expressions of the form $\sum_{i=1}^n\pm r_i$.  
The space $\M_n$ dominates all the spaces $\M_r$, but each $\M_r$ has a dense open 
subset which is smooth and diffeomorphic to a smooth dense open subset of $\M_n$. 
We do not know whether or not the real analytic space $\M_n$ or the Chow quotient 
$Y^n//G$ is smooth.  As in \cite{Hu}, a point in the space $\M_n$ can be represented 
as a ``bubble'' polygon.    

	Another interesting observation, which emphasizes the importance of the 
involution $\t$ on the complex flag manifolds, is that, with just a little effort, 
one can recover the Schubert calculus for the quaternionic Grassmannians from 
that of the complex ones, restricting the attention to the cycles preserved by 
$\t$ and considering their $\t$-fixed point subsets.  Previously, the quaternionic 
Schubert calculus was considered in \cite{PragaczR} and was dealt with using 
different methods.  
\bigskip

\section{Generalized action-angle coordinates on $\Mr$}

In this section we describe certain local ``action-angle coordinates" in a open dense subset of $\Mr$, which are analogous to those defined in \cite{KM} for the $\3$ case.

	Let $\ell_i, \ i=3,...,n-1$, denote the length of each of the $n-3$ diagonals $d_i$ of $[P]$ in $\Mr$, as defined in Section $4$.  We will call $[P]\in \Mr$ {\it generic} if $\ell_i\neq 0$ and $\ell_i+r_{i+2}\neq \ell_{i+1}$, for all $i=3,...,n-2$.

	For every generic $[P]\in\Mr$ we construct a canonical planar $n$-gon $[P_c]\in\Mr$ as follows.  Choose any representative $P$ of $[P]$ having its first vertex at the origin of $\R^5$, and let $\Pi^P$ be the $2$-dimensional subspace of $\R^5$ spanned by the first three vertices of $P$.  Then one may use the copy of $SO(4)$ fixing $\Pi^P$ to rigidly move the fourth vertex of $P$ to $\Pi^P$ in such a way that the line segment joining $v_2$ and $v_4$ intersects the $1$-dimensional subspace of $\R^5$ containing $d_3$, the first diagonal of $P$.  The generic character of $P$ ensures that one may repeat this procedure enough times so as to eventually obtain a unique planar polygon $P_c$ lying entirely on $\Pi^P$, and such that the segment connecting $v_i$ and $v_{i+2}$ always intersects the line in $\R^5$ containing $d_{i+1}$.  Clearly, the correspondence $[P]\rightarrow[P_c]$ is well-defined, and $[P_c]$ is unique for generic $[P]$.

	Let $L_{\bf \ell}\subseteq\Mr$ be a level set of the $n-3$ length functions $\ell_i$.  We now describe $3n-12$ local ``angle" coordinates for $L_{\bf \ell}$ within the open subset of $\Mr$ consisting of generic $n$-gons.  Let $[P_c]$ denote the canonical planar $n$-gon associated to some generic $[P]\in L_{\bf \ell}$.  Assume that $P_c\subseteq\Pi^P$ is chosen so that $v_1$ lies at the origin of $\R^5$ and let $K_i$ denote the copy of $SO(4)\subseteq SO(5)$ fixing $d_i, \ i=3,...,n-1$.  Consider all possible deformations of $P_c$ obtained by moving its second vertex, $v_2$, by $K_3$ and fixing its remaining vertices.  Let $F_2$ denote the subgroup of $K_3$ fixing $v_2$ (and hence all of $P_c$).  As $F_2$ is isomorphic to $SO(3)$, we see that we may parameterize the above family of deformations of $P_c$ in $\R^5$ by $K_3/F_2\simeq SO(4)/SO(3)\simeq S^3$.  (Indeed, whatever the deformation of $P_c$, $v_2$ is constrained to the intersection of two $4$-spheres of radius $r_1$ and $r_2$ respectively).  Note however, that the collection of $\R^5$-polygons obtained by deforming $P_c$ as above is not in one-to-one correspondence with generic $[P]$ in $\Mr$;  that is, there are polygons in this collection which differ by a rigid motion of $\R^5$.  Given the nature of the deformations of $P_c$ being considered, any such rigid motion must fix the plane $\Pi^{P_c}$.  Consequently, deformations of $[P_c]$ in $\Mr$ by bending along its first diagonal are parameterized by $SO(3)\setminus SO(4)/SO(3)\simeq S^3/SO(3)\simeq [0,\pi]$, and so we get a first angle coordinate for $L_{\bf \ell}$.
  
	Let $P_c^3$ be a generic deformation of $P_c$ obtained from $K_3$ by bending along $d_3$ as described above, and consider all possible deformations of $P_c^3$ obtained from $K_4$ by bending along $d_4$.  Since the first four vertices of $P_c^3$ span a copy of $\R^3$ fixed by some $F_3\simeq SO(2)\subseteq K_4$, we see that these deformations of $P_c^3$ may be parameterized by $SO(4)/SO(2)$.  As above, corresponding deformations of $[P_c^3]$ are then parameterized by $SO(3)\setminus SO(4)/SO(2)\simeq \overline{D^2}$, where $SO(3)$ fixes $\Pi^{P_c}$ and $\overline{D^2}$ denotes a closed disk.  Thus we obtain two more (independent) ``generalized" angle variables.  

	Finally, let $P_c^i$ denote a generic deformation of $P_c^{i-1}$ obtained from $K_i$ by bending along $d_i$, $i=4,...,n-1$.  Since the first $i+1$ vertices of $P_c^i$ span at least $\R^4$, we need all of $K_{i+1}\simeq SO(4)$ to parameterize all possible rigid deformations of $P_c^i$ by bending about $d_{i+1}$.  Corresponding deformations of $[P_c^i], \ i=4,...,n-1$, are hence parameterized by $SO(3)\setminus SO(4)$, yielding $3(n-5)$ new generalized angle coordinates.  Hence, we obtain a total of $3n-12$ generalized angle variables, which together with the $n-3$ action variables prescribed by the lengths $\ell_i$, yield a set of local generalized coordinates for the dense open subset of $\Mr$ consisting of generic polygons. 

	We remark that it is easy to see from the above description how to stratify $L_{\ell}$ as a disjoint union of smooth manifolds.

\bibliographystyle{alpha}
\bibliography{FLQuatPol}

\vskip 0.2in

\noindent
Department of Mathematics \\ 
University of Arizona \\
Tucson, AZ 85721-0089 
\medskip

\noindent
{\tt foth@math.arizona.edu} \\
{\tt lozano@math.arizona.edu}
\vskip+.5in

\end{document}